[1994/12/01]
\documentclass[12pt, reqno]{amsart}
\usepackage{a4wide}
\usepackage[english, activeacute]{babel}
\usepackage{amsmath,amsthm,amsxtra}
\usepackage{epsfig}
\usepackage{amssymb}
\usepackage{latexsym}
\usepackage{amsfonts}
\usepackage[colorlinks=true]{hyperref}
\hypersetup{linkcolor=red,citecolor=blue,filecolor=dullmagenta,urlcolor=red} 

\title{On the asymptotic behavior of solutions to the Benjamin-Ono equation}
\author{Claudio Mu\~noz}

\author{Gustavo Ponce}
\thanks{}
\address{CNRS and Departamento de Ingenier\'ia Matem\'atica DIM-CMM UMI 2807-CNRS \\ Universidad de Chile, Santiago, Chile}
\email{cmunoz@dim.uchile.cl}
\address{Department of Mathematics \\ University of California-Santa Barbara, CA 93106-USA}
\email{ponce@math.ucsb.edu}
\date{\today}

\chardef\bslash=`\\ 

\newtheorem{thm}{Theorem}[section]

\newtheorem{lem}[thm]{Lemma}

\theoremstyle{remark}
\newtheorem{rem}{Remark}[section]

\numberwithin{equation}{section}

\newcommand{\R}{\mathbb{R}}

\newcommand{\Z}{\mathbb{Z}}

\newcommand{\sgn}{\operatorname{sgn}}

\newcommand{\be}{\begin{equation}}
\newcommand{\ee}{\end{equation}}
\newcommand{\bp}{\begin{proof}}
\newcommand{\ep}{\end{proof}}
\newcommand{\bel}{\begin{equation}\label}
\newcommand{\eeq}{\end{equation}}
\newcommand{\bea}{\begin{eqnarray}}
\newcommand{\eea}{\end{eqnarray}}
\newcommand{\bee}{\begin{eqnarray*}}
\newcommand{\eee}{\end{eqnarray*}}
\newcommand{\ben}{\begin{enumerate}}
\newcommand{\een}{\end{enumerate}}

\newcommand{\eval}[2][\right]{\relax
  \ifx#1\right\relax \left.\fi#2#1\rvert}

\begin{document}
\begin{abstract}

We prove that the limit infimum, as time $\,t\,$  goes to infinity, of any uniformly bounded in time $H^1\cap L^1$ solution to the Benjamin-Ono equation converge to zero locally in an increasing-in-time region of space of order $\,t/\log t$. 
Also for a solution with a mild $L^1$-norm growth in time,  its limit infimum must converge to zero, as time goes to infinity, locally in an increasing on time region of space of order  depending of the rate of growth of its $L^1$-norm.  In particular, we discard the existence of breathers and other  solutions for the BO model moving with a speed \lq\lq slower"  than a soliton.
\end{abstract}
\maketitle \markboth{BO decay estimates}{C. Mu\~noz-G. Ponce}
\renewcommand{\sectionmark}[1]{}

\section{Introduction and main results}


This work is concerned with the  initial value problem (IVP) associated to Benjamin-Ono equation (BO) 
\be\label{BO} 
\begin{aligned}
\begin{cases}
& \!\! \partial_t u + \partial_{x}(\partial_x\mathcal H u +u^2) = 0,\qquad (t,x) \in \R\times \R,\\
&\!\! u(x,0)=u_0(x),
 \end{cases}
\end{aligned}
\ee
where $\mathcal H$ denotes the Hilbert transform
\begin{equation}
\label{hita}
\begin{aligned}
\mathcal H f(x) :=&~ \frac{1}{\pi} {\rm p.v.}\Big(\frac{1}{x}\ast f\Big)(x)
\\
:=& ~ \frac{1}{\pi}\lim_{\epsilon\downarrow 0}\int\limits_{|y|\ge \epsilon} \frac{f(x-y)}{y}\,dy=(-i\,\sgn(\xi) \widehat{f}(\xi))^{\vee}(x).
\end{aligned}
\end{equation}

The BO equation was first deduced by Benjamin
\cite{Be} and Ono \cite{On} as  a model for long internal gravity
waves in deep stratified fluids. Later, it was also shown to be   a
completely integrable system (see \cite{AS} and
references therein). Thus,    the BO equation has followed the same historical pattern of the celebrated 
Korteweg-de Vries (KdV) equation \cite{KdV} 
\begin{equation}\label{kdv}
\partial_t u+\partial_x^3u+u\,\partial_x u=0.
\end{equation}
 As a  completely integrable system solutions of the IVP \eqref{BO} satisfy infinitely many conservation laws. The first three of them are
\be
\label{cl}
\begin{aligned}
&I_1(t)=\int_{-\infty}^{\infty} u(x,t)dx =I_1(0),\;\;\;\;\;\;\;\;\;I_2(t)=\int_{-\infty}^{\infty} u^2(x,t)dx =I_2(0),\\
\\
&E(t)=\int_{-\infty}^{\infty} \Big( (D^{1/2}u)^2-\frac{1}{3}\,u^3\Big)(x,t)dx=E(0), \qquad D^s:=\mathcal{F}^{-1} |\xi|^s \mathcal F.
\end{aligned}
\ee

The local well-posedness (LWP) and global well-posedness (GWP) of the IVP for the BO equation in classical Sobolev spaces $\,H^s(\R)= \left(1-\partial^2_x\right)^{-s/2} L^2(\R),\,s\in\R\,$ has been extensively studied: in \cite{ABFS} and  \cite{Io1} LWP was obtained for $s>3/2$, in \cite{Po} GWP was proven for $s\geq 3/2$, LWP was established in \cite{KoTz1} for $s>5/4$ and in  \cite{KeKo} for $s>9/8$, in \cite{Ta} GWP was demonstrated  for $s\geq 1$, in \cite{BuPl} LWP was found to hold for $s>1/4$ and in
\cite{IoKe} GWP was presented in $H^0(\R)=L^2(\R)$, (for further details  and results regarding the well-posedness of the IVP associated to the BO equation \eqref{BO} in $H^s(\R)$  see \cite{MoPi}, \cite{LiPo}).
More recently, in \cite {IT} a different proof of the result in \cite{IoKe} was provided as well as some description on the long time behavior of appropriate small solution to \eqref{BO}.

\medskip

The above results have been obtained by \lq\lq compactness methods", since one can only show  that the map \emph{data $\to$ solution} is just locally continuous. In fact, in \cite{MoSaTz} it was established that for any $s\in\R$ the flow map $u_0\to u$ from $H^s(\R)$ to $\,C([-T,T]:H^s(\R))$ is not locally of class $C^2$  (see also \cite{KoTz2}).

\medskip

It was first deduced in \cite{Io1} and \cite{Io2} that, in
general,  decay of polynomial type is not preserved by the solution flow of the BO equation. The results in \cite{Io1} and \cite{Io2} have been extended to fractional order
weighted Sobolev spaces and have shown to be optimal  in \cite{FP} and \cite{FLP}. 

\medskip

To make some of the above statements precise, we  introduce the weighted Sobolev spaces
\begin{equation}
\label{spaceZ}
Z_{s,r} := H^s(\R) \cap L^2(|x|^{2r}dx),\,\,\, s, r \in \R ,
\end{equation}
and
\begin{equation}
\label{spaceZdot}
\dot Z_{s,r} :=\{ f\in H^s(\R)\cap L^2(|x|^{2r}dx)\,:\,\widehat {f}(0)=0\},\;\;\;\;\;\;\;s,\,r\in\R.
\end{equation}

Thus, in particular:

\medskip

{\it(i)}\,\cite{FP} \,Let $u\in C([-T,T] : Z_{2,2})$ be a solution of the IVP for the equation \eqref{BO}. If   there exist  two different times
$\,t_1, t_2\in [-T,T]$ such that
\begin{equation}
\label{2timesw}
u(\cdot,t_j)\in Z_{5/2,5/2},\;j=1,2,\;\text{then}\;\;\widehat {u}_0(0)=0,\;\,(\text{so}\;\, u(\cdot, t)\in  \dot Z_{5/2,5/2}).
\end{equation}
The solution flow of \eqref{BO} preserves the class $\,\dot Z_{s,r}$ for $s\geq r\in (1/2,7/2)$.

\vskip.05in

{\it(ii)}\,\cite {FLP} \, Let $u\in C([-T,T] : \dot Z_{3,3})$ be a solution of the IVP for the equation  \eqref{BO}. If   there exist  three different times
$\,t_1, t_2, t_3\in [-T,T]$ such that
\begin{equation}
\label{3timeus}
u(\cdot,t_j)\in \dot Z_{7/2,7/2},\;\;j=1,2,3,\;\;\text{then}\;\;\;u(x,t)\equiv 0. \end{equation}

\vskip.05in

In this work we are concerned with the long time behavior of solutions of the IVP \eqref{BO}. More precisely, we are interested in the \lq\lq location" of the $\,H^{1/2}$-norm of the solution which is globally bounded, due to the conservation laws  \eqref{cl},  as time evolves. 

\medskip

Here we shall asume that:  $\exists \;a\in [0,1/2) \;\;\,\exists\; c>0\;$ such that $\,\,\forall \,T>0$
\be\label{L1}
\sup_{t\in[0,T]} \int_{-\infty}^{\infty} |u(x,t)|dx\leq c_0\,\langle T\rangle ^a, \qquad\;\;\;\;\;\; \langle T\rangle:= (1+T^2)^{1/2}.
\ee

From  the original proof of the so called local smoothing effect found in \cite{Ka}, see also \cite{KF},  several modified versions of this weighted energy estimates argument have been quite successfully used, see for example the works \cite {MM} and \cite{MM1}. These works are centered in the understanding of the long time behavior around stable and unstable generalized KdV solitons, i.e., positive speed solutions.

\medskip

More related with the problem addressed here we have the variants found in the recent works \cite{KMM,KMM1,KMPP}. However, in the case treated here even under the strongest hypothesis $\,a=0$ in \eqref{L1} we need a weight outside the cut-off function $\,\phi(\cdot)$, see \eqref{weight1}. This is a consequence of the dispersive relation of the BO which  is not strong enough to apply the process as in \cite{KMM1} or \cite{MuPo}. This extra-weight allows us to close the estimate in a weaker form than those obtained in \cite{KMM1} and \cite{MuPo}, since it involves the lim inf instead of the standard lim. 

\medskip

 Also, to implement our approach  we have to rely on estimates obtained in \cite{KeMa}, which requires a very specific cut-off function 
$\,\phi$, see Lemma \ref{KM} and Lemma \ref{KM11}. In fact, one observes that $\,\phi'(x)$ is a multiple of the soliton solution of the BO equation, see \eqref{soliton}.

\medskip

Our main result in this work is the following theorem:

\begin{thm}\label{TH2}
Let $u=u(x,t)$ be  a solution to the IVP  \eqref{BO} such that
 \begin{equation}
 \label{m11}
 u\in C(\R:H^{1}(\R)) \cap L^{\infty}_{loc}(\R:L^1(\R))
\end{equation}
satisfying \eqref{L1}. Then,
\be\label{main}
\int_{10}^{\infty}\; \frac{1}{ t\,\log t}\;\Big(\int_{-\infty}^{\infty}\phi'\Big(\frac{x}{\lambda(t)}\Big)(u^2+(D^{1/2}u)^2)(x,t)dx\Big)\,dt<\infty.
\ee
Hence, 
\be
\label{abc2}
\liminf_{t\uparrow \infty}\; \int_{-\infty}^{\infty} \,\phi'\Big(\frac{x}{\lambda(t_n)}\Big)(u^2+(D^{1/2}u)^2)(x,t)\,dx= 0,
\ee
with
\be\label{def1}
\lambda(t)=\frac{c\,t^b}{\log t },\;\;\;\;\;\;\;a+b=1,\;\;\;\;\;\text{and}\; \;\;\;\phi'(x)=\frac{1}{1+x^2},
\ee
for any fixed $\,c>0$.
\end{thm}

\begin{rem}
Theorem \ref{TH2}  tells us that there exists a sequence of times $\{t_n\,:\,n\in\Z^+\}$ with $t_n\uparrow \infty$ as $n\uparrow \infty$ such that
\be\label{asymp22}
\lim_{n\uparrow \infty} \,\int_{|x|\leq \frac{c\,t^{b}_n}{\log(t_n)}}\;
(u^2+(D^{1/2}u)^2)(x,t_n)dx=0.
\ee

Since the BO equation is time reversible the same result holds for the time interval $\,(-\infty,-10]$. We do not know how to prove
\eqref{abc2} with $\,\lim\,$ instead of $\,\lim\inf$. 
Notice that in the case $\,a=0$ in \eqref{L1} one has $\,b=1$. Thus, roughly speaking, in this case such result would comply with  the statement of the so-called soliton resolution conjecture.

\end{rem}

\begin{rem}
The result in Theorem \ref{TH2} shows that the BO can not possess time periodic non-trivial solutions (in particular breathers). However, from the argument in \cite{MuPo} this result will follow by observing that a solution with mild space decay  satisfies
\[
\frac{d\;}{dt}\int\,xu(x,t)dx=\int\,u^2(x,t)dx.
\]
\end{rem}

\begin{rem} The assumption \eqref{m11} can be deduced by assuming that $u_0\in Z_{s,r}$ with $\,s\geq r>1/2$, see \cite{FP}.

\end{rem}

\begin{rem}
Note that solitons 
\be
\label{soliton}
u(x,t)=Q_c(x-ct),\;\;\;\;\;\;\;\;\;\;Q_c(x)=\frac{4c}{1+c^2x^2},
\ee are in the class \eqref{m11}, and they also satisfy \eqref{abc2}.
\end{rem}

\begin{rem}
Note that Theorem \ref{TH2} works for nonintegrable perturbations of the BO equation, as long as the solution is sufficiently small.
\end{rem}

The rest of this paper is organized as follows: Section 2 contains all the preliminary results needed in the proof of Theorem \ref{TH2} which will be given in Section 3.

\bigskip

\section{Preliminary estimates }

We shall take 
\be
\label{def11}
\phi(x) = \frac{\pi}{2} +\tan^{-1}(x),\;\;\;\; \;\;\;\;\,\lambda(t)=\frac{c\,t^{b}}{\log t },\;\;\;c>0,\;\;\;0<b= 1-a,
\ee
and recall that
\be\label{Hilbert}
\phi'(x) =\frac{1}{1+x^2},\;\;\;\;\;\;\mathcal H \phi'(x)= \frac{x}{1+x^2},\;\;\;\;\;\mathcal H \phi''(x)= \frac{1-x^2}{(1+x^2)^2},
\ee
(to simplify the argument without loss of generality we shall assume $c=1$ in \eqref{def11}).

\medskip

In the proof of Theorem \ref{TH2}  we shall use the following inequalities :

\begin{lem} \label{KM}
For any $f\in H^1(\R)$, 
\be
\label{KM1}
\int (\mathcal H\partial_x f) \,f\, \phi'\Big(\frac{x}{\lambda(t)}\Big)dx \leq \frac{c}{\lambda(t)}\,\int \,f^2\,\phi'\Big(\frac{x}{\lambda(t)}\Big)dx,
\ee
\be\label{KM2}
\Big |\,\int (\mathcal H\partial_x f) \,(\partial_x f)\,\phi\Big(\frac{x}{\lambda(t)}\Big)dx\Big|\leq \frac{c}{\lambda(t)}\,\int \,f^2\,\phi'\Big(\frac{x}{\lambda(t)}\Big)dx,
\ee
and for $u\in L^2$,  $ \varphi$ weight,
\be\label{comm}
\| D^{1/2}\big[ D^{1/2};\varphi\big] u\|_2\leq c \| \widehat{\partial_x \varphi}\|_1\,\|u\|_2.
\ee
\end{lem}

\begin{lem} \label{KM11}

If $\,u\in C(\R:H^{1/2}(\R))$ is a strong solution of the BO equation in \eqref{BO}, then
\be\label{key}
\int_{-\infty}^{\infty} |u(x,t)|^3\,\phi'\Big(\frac{x}{\lambda(t)}\Big)dx\leq c\,\int_{-\infty}^{\infty} (u(x,t))^2\,\phi'\Big(\frac{x}{\lambda(t)}\Big)dx,
\ee
where $\,c=c(\|u_0\|_2\,;\,\| D^{1/2}u_0\|_2)$.
\end{lem}

The estimates \eqref{KM1} and \eqref{KM2} were proven in \cite{KeMa} (Lemma 2 and Lemma 3, resp.). For the proof of \eqref{comm} 
we refer to a similar one in  \cite{DMP}.  For the proof of \eqref{key} we refer to \cite{KeMa}  (Appendix 2).

\section{Proof of Theorem \ref{TH2}}

Multiplying the equation \eqref{BO} by 
\be\label{weight1}
\frac{1}{ t^a\,\log^2(t)}\,\phi\Big(\frac{x}{\lambda(t)}\Big),
\ee
one gets that
\begin{equation}
\label{eq1}
\begin{aligned}
& \frac{d}{dt} \int\frac{1}{ t^a\,\log^2(t)}\,\phi\Big(\frac{x}{\lambda(t)}\Big)\,u(x,t)dx \\
\\
&\underbrace{-\int \Big(\frac{1}{ t^a\,\log^2(t)}\Big)'\,\phi\Big(\frac{x}{\lambda(t)}\Big)\,u(x,t) dx}_{=:A_1(t)}  \\
\\
&+
\underbrace{\int \frac{1}{ t^a\,\log^2(t)}\,\frac{x}{\lambda(t)}\,\frac{\lambda'(t)}{\lambda(t)}\phi'\Big(\frac{x}{\lambda(t)}\Big)\,u(x,t) dx}_{=:A_2(t)}
\\
\\
& -
\underbrace{\int \frac{1}{ t^a\,\log^2(t)}\,\phi\Big(\frac{x}{\lambda(t)}\Big)\,\partial_x^2\,\mathcal H u(x,t)\,dx}_{=:A_3(t) }\\
\\
&-\underbrace{\int \frac{1}{ t^a\,\log^2(t)}\,\phi\Big(\frac{x}{\lambda(t)}\Big)\,\partial_x(u^2(x,t))\,dx}_{=:A_4(t) }=0.
\end{aligned}
\end{equation}

\vskip.05in

Observing that
$$
\Big(\frac{1}{ t^a\,\log^2(t)}\Big)'\sim \ \frac{1}{ t^{1+a}\,\log^2(t)},\;\;\;\;\;\;\;\;\text{as}\;\;\;\;t\uparrow \infty,
$$
from the assumption \eqref{L1} it follows that
\be\label{a1}
A_1(t)\sim \frac{1}{t\,\log^2(t)}\in L^1([10,\infty)).
\ee

\vskip.05in 

Since 
\be\label{a2}
\frac{\lambda'(t)}{\lambda(t)}\sim\frac{1}{t},\;\;\;\;\;\;\;\;\text{as}\;\;\;\;t\uparrow \infty,
\ee
and
\be\label{aaa2}
\frac{x}{\lambda(t)}\,\phi'\Big(\frac{x}{\lambda(t)}\Big)\in L^{\infty}(\R)\;\;\text{uniformly in} \;t\in[10,\infty),
\ee
one sees that 
\be\label{a3}
A_2(t)\sim \frac{1}{t\,\log^2(t)}\in L^1([10,\infty)).
\ee
Similarly, after integration by parts, using \eqref{Hilbert} it follows that
\be
\label{abcd}
\partial_x^2\mathcal H \phi\Big(\frac{x}{\lambda(t)}\Big)=\frac{1}{\lambda^2(t)}\,\frac{1-(x/\lambda(t))^2}{(1+(x/\lambda(t))^2)^2},
\ee
with
$$
\frac{1-x^2}{(1+x^2)^2}\in L^{\infty}(\R).
$$
Thus,
\be\label{a4}
A_3(t)\sim \frac{1}{(\lambda(t))^2\,\log^2(t)}\in L^1([10,\infty)),\;\;\;\;\;\text{if}\;\;\;\;b>1/2.
\ee
Hence, integrating in time the identity \eqref{eq1}, and combining the above results and some integration by parts one concludes that
\be\label{a5}
\begin{aligned}
&\int_{10}^{\infty} \,\int_{-\infty}^{\infty} \frac{1}{ t^a\,\log^2(t)}\,\frac{1}{\lambda(t)}\,\phi'\Big(\frac{x}{\lambda(t)}\Big)u^2(x,t)dx\,dt\\
\\
&=\int_{10}^{\infty} \frac{1}{ t\,\log(t)}\;\Big(\int_{-\infty}^{\infty}\phi'\Big(\frac{x}{\lambda(t)}\Big)u^2(x,t)dx\Big)\,dt<\infty.
\end{aligned}
\ee
\vskip.05in

Next, we multiplying the equation \eqref{BO} by 
\be\label{weight}
\frac{1}{ t^a\,\log^2(t)}\,\phi\Big(\frac{x}{\lambda(t)}\Big) \,u(x,t)
\ee
to get, after some integration by parts, that
\begin{equation}
\label{111}
\begin{aligned}
& \frac{1}{2}\,\frac{d}{dt} \int\frac{1}{ t^a\,\log^2(t)}\,\phi\Big(\frac{x}{\lambda(t)}\Big)\,u^2(x,t)dx \\
\\
&\underbrace{- \frac{1}{2}\int \Big(\frac{1}{ t^a\,\log^2(t)}\Big)'\,\phi\Big(\frac{x}{\lambda(t)}\Big)\,u^2(x,t) dx}_{=:B_1(t)}  \\
\\& +
\frac{1}{2}\underbrace{\int \frac{1}{ t^a\,\log^2(t)}\,\frac{x}{\lambda(t)}\,\frac{\lambda'(t)}{\lambda(t)}\phi'\Big(\frac{x}{\lambda(t)}\Big)\,u^2(x,t) \,dx}_{=:B_2(t)}
\\
\\
& 
-\underbrace{\int\frac{1}{ t^a\,\log^2(t)}\,\phi\Big(\frac{x}{\lambda(t)}\Big)\,(\mathcal H \partial^2_xu)(x,t)\,u(x,t)\,dx}_{=:B_3(t)}
\\
\\
&+\frac{2}{3}\underbrace{\int \frac{1}{ t^a\,\lambda(t)\,\log^2(t)}\,\phi'\Big(\frac{x}{\lambda(t)}\Big)\,u^3(x,t)\,dx}_{=:B_4(t) }=0.
\end{aligned}
\end{equation}

As before one has that
$$
\Big(\frac{1}{ t^a\,\log^2(t)}\Big)'\sim \ \frac{1}{ t^{1+a}\,\log^2(t)},\;\;\;\;\;\;\;\;\text{as}\;\;\;\;t\uparrow \infty,
$$
so from the preservation of the $L^2$-norm of the solution it follows that
\be\label{a11}
B_1(t)\sim \frac{1}{t^{1+a}\,\log^2(t)}\in L^1([10,\infty)),\;\;\;\;\;\forall a\geq0.
\ee
Similarly,
\be\label{a12}
B_2(t)\sim\frac{\lambda'(t)}{t^a\,\log^2(t)\,\lambda(t)}\sim  \frac{1}{t^{1+a}\,\log^2(t)}\in L^1([10,\infty)),\;\;\;\;\;\forall a\geq 0.
\ee

From Lemma \ref{KM11} one finds that
\be
\label{new1}
B_4(t)\leq \frac{c}{ t\,\log(t)}\;\Big(\int_{-\infty}^{\infty}\phi'\Big(\frac{x}{\lambda(t)}\Big)u^2(x,t)dx\Big),
\ee
which by \eqref{a5} one has
\be\label{aa11}
B_4(t)\in L^1([10,\infty)),\;\;\;\;\;\forall a\geq0.
\ee

It remains to consider $\,B_3(t)\,$ in \eqref{111}
 where
 \be
 \label{p1}
 B_3(t)=\,-\,\frac{1}{t^a\,\log^2(t)}\,D_3(t),
 \ee
with
\be\label{p2}
\begin{aligned}
D_3(t)&=\int \mathcal H\partial_xu\,\partial_xu\,\phi\Big(\frac{x}{\lambda(t)}\Big) dx +
 \frac{1}{\lambda(t)}\,\int \mathcal H\partial_x u\,u\,\phi'\Big(\frac{x}{\lambda(t)}\Big) dx\\
 \\
 & =D_{3,1}(t)+\frac{D_{3,2}(t)}{\lambda(t)}.
 \end{aligned}
 \ee

By Lemma \ref{KM} (estimate \eqref{KM2})
\be\label{p3}
|D_{3,1}(t)|= \Big| \int\,\mathcal H\partial_x u\,\partial_xu \,\phi\Big(\frac{x}{\lambda(t)}\Big)dx\,\Big|\leq \frac{1}{\lambda(t)}\int u^2\,
\phi'\Big(\frac{x}{\lambda(t)}\Big)dx,
\ee
with
$$
\phi'\Big(\frac{x}{\lambda(t)}\Big)=\frac{1}{\lambda(t)}\,\frac{1}{1+(x/\lambda(t))^2}.
$$
Hence,
\be\label{p5}
|D_{3,1}(t)|\leq \frac{1}{\lambda^2(t)}\|u_0\|_2^2,
\ee
so
\be\label{p7}
\frac{1}{t^a\,\log^2(t)} |D_{3,1}(t)|\sim \frac{1}{t^a\,\log^2(t)\,\lambda^2(t)}\in L^1([10,\infty)),\;\;\;\;\;\text{since}\;\;\;\;\;a+2b>1.
\ee

Also writing 
\be
\label{p7a}
\begin{aligned}
D_{3,2}(t)&=\int Du\,u\,\phi'\Big(\frac{x}{\lambda(t)}\Big)dx\\
&=\int D^{1/2}u\,D^{1/2}u\,\phi'\Big(\frac{x}{\lambda(t)}\Big)dx+ \int D^{1/2}u\,\Big[D^{1/2};\phi'\Big(\frac{x}{\lambda(t)}\Big)\Big] u\,dx\\
&=D_{3,2,1}(t)+D_{3,2,2}(t).
\end{aligned}
\ee

The total contribution of the term involving $ D_{3,2,1}(t)$ is
\be
\label{007}
\frac{1}{t^a\,\log^2(t)\,\lambda(t)}\,\int (D^{1/2}u)^2\,\phi'\Big(\frac{x}{\lambda(t)}\Big)\,dx.
\ee

From Lemma \ref{KM} (estimate \eqref{comm}) it follows that
\be\label{p9}
\begin{aligned}
|D_{3,2,2}(t)|=\Big|\int u\,D^{1/2}\,\Big[D^{1/2};\phi'\Big(\frac{x}{\lambda(t)}\Big)\Big] u\,dx \Big|\leq c\,\|\widehat{\partial_x\phi'\Big(\frac{\cdot}{\lambda(t)}\Big)}\|_1 \|u\|_2.
\end{aligned}
\ee
Since
\[
\widehat{\partial_x\phi'\Big(\frac{\cdot}{\lambda(t)}\Big)}(\xi) =\frac{1}{\lambda(t)}\,\widehat {\phi''\Big(\frac{\cdot}{\lambda(t)}\Big)}(\xi)
\]
with 
\[
\,\phi''(x)=-\frac{2x}{(1+x^2)^2},
\]
using that in general
\[
\|\,\widehat{f}\|_1\leq c\,(\|f\|_2 +\|\,f'\|_2),
\]
 one has that
\[
\Big\| \widehat{\partial_x\phi'\Big(\frac{\cdot}{\lambda(t)}\Big)}(\xi) \Big\|_1\sim \frac{1}{(\lambda(t))^{1/2}}\;\;\;\;\;\;\;\;\text{as}\;\;\;\;t\uparrow \infty.
\]
So the total contribution of the term $\,D_{3,2,2}(t)$ in \eqref{p1} and \eqref{p2} is bounded
by 
\be\label{p7b} \frac{1}{t^a\,\log^2(t)\,\lambda^{3/2}(t)}\in L^1([10,\infty)),\;\;\;\;\;\text{since}\;\;\;\;\;a+\frac{3}{2}b>1.
\ee

Gathering the above information we can conclude that 
\be\label{aa5}
\begin{aligned}
&\int_{10}^{\infty} \,\int_{-\infty}^{\infty} \;\frac{1}{ t^a\,\log^2(t)}\,\frac{1}{\lambda(t)}\,\phi'\Big(\frac{x}{\lambda(t)}\Big)(D^{1/2}u)^2(x,t)dx\,dt\\
\\
&=\int_{10}^{\infty}\; \frac{1}{ t^{a+b}\,\log(t)}\;\Big(\int_{-\infty}^{\infty}\phi'\Big(\frac{x}{\lambda(t)}\Big)(D^{1/2}u)^2(x,t)dx\Big)\,dt<\infty.
\end{aligned}
\ee
Therefore, collecting the information in \eqref{a5} and \eqref{aa5}
\be\label{aaa5}
\begin{aligned}
&\int_{10}^{\infty} \,\int_{-\infty}^{\infty} \;\frac{1}{ t^a\,\log^2(t)}\,\frac{1}{\lambda(t)}\,\phi'\Big(\frac{x}{\lambda(t)}\Big)(u^2+(D^{1/2}u)^2)(x,t)dx\,dt\\
\\
&=\int_{10}^{\infty}\; \frac{1}{ t^{a+b}\,\log(t)}\;\Big(\int_{-\infty}^{\infty}\phi'\Big(\frac{x}{\lambda(t)}\Big)(u^2+(D^{1/2}u)^2)(x,t)dx\Big)\,dt<\infty,
\end{aligned}
\ee
with 
\[
\,\phi'(x)=\frac{1}{1+x^2}.
\]
 Since $\,a +b=1$, and 
\be
\label{nonint}
\eta(t):=\frac{1}{t\,\log(t)}\notin  L^1([10,\infty)),
\ee
from \eqref{a5} it follows that  there exists $\,\{t_n\;:\;n\in\mathbb N\}$ with $\,t_n\uparrow \infty$ such that
\be
\label{abc1}
F(t_n)=:\int_{-\infty}^{\infty} \,\phi'\Big(\frac{x}{\lambda(t_n)}\Big)(u^2+(D^{1/2}u)^2)(x,t_n)dx\to 0\;\;\;\;\;\text{as}\;\;\;\;n\uparrow \infty.
\ee
Moreover, one can see that for any $\,k\in\mathbb N$ there exists $\,t_k\in[2^k, 2^{k+1})$ such that

\be \label{est1}
F(t_k)\leq \frac{1}{\log(k)}.
\ee

Above we have proved that
\be
\label{abc22}
 \int_{-\infty}^{\infty} \,\phi'\Big(\frac{x}{\lambda(t_n)}\Big)(u^2+(D^{1/2}u)^2)(x,t_n)dx\to 0\;\;\;\;\;\text{as}\;\;\;\;n\uparrow \infty
\ee
with
\be\label{def22}
\lambda(t)=\frac{c\,t^b}{\log(t)},\;\;\;\;\;\text{and}\;\;\;\;\;\;a+b=1.
\ee
In particular, one has that
\be\label{asymp222}
\lim_{n\uparrow \infty} \,\int_{|x|\leq \frac{c\,t^{b}_n}{\log(t_n)}}\;
(u^2+(D^{1/2}u)^2)(x,t_n)dx=0
\ee
which yields the desired result.


\end{document}